\documentstyle[bezier]{article}

\makeatletter
\def\eqnarray{%
  \stepcounter{equation}%
  \let\@currentlabel=\theequation
  \global\@eqnswtrue
  \global\@eqcnt\z@
  \tabskip\@centering
  \let\\=\@eqncr
  $$\halign to \displaywidth\bgroup\@eqnsel\hskip\@centering
  $\displaystyle\tabskip\z@{##}$&\global\@eqcnt\@ne
  \hfil$\displaystyle{{}##{}}$\hfil
  &\global\@eqcnt\tw@$\displaystyle\tabskip\z@{##}$\hfil
  \tabskip\@centering&\llap{##}\tabskip\z@\cr}
\makeatother

\makeatletter
  \renewcommand{\theequation}{%
        \thesection.\arabic{equation}}
  \@addtoreset{equation}{section}
\makeatother

\begin{document}

\newtheorem{th}{Donotwrite}[section]

\newtheorem{definition}[th]{Definition}
\newtheorem{theorem}[th]{Theorem}
\newtheorem{proposition}[th]{Proposition}
\newtheorem{lemma}[th]{Lemma}
\newtheorem{corollary}[th]{Corollary}
\newtheorem{remark}[th]{Remark}
\newtheorem{example}[th]{Example}

\newfont{\germ}{eufm10}

\def\Aff{\mbox{\sl Aff}\,}
\def\btilde{\tilde{b}}
\def\cd{\cdots}
\def\eps{\epsilon}
\def\goth#1{\mbox{\germ #1}}
\def\ot{\otimes}
\def\Proof{\noindent{\sl Proof.}\quad}
\def\Q{{\bf Q}}
\def\qed{~\rule{1mm}{2.5mm}}
\def\R{{\cal R}}
\def\slchap{\widehat{\goth{sl}}_n}
\def\sln{\goth{sl}_n}
\def\sln-{\goth{sl}_{n-1}}
\def\ub(#1,#2){\underbrace{#1\ot\cd\ot#1}_{#2}}
\def\ubb(#1,#2){\underbrace{#1\cd#1}_{#2}}
\def\veps{\varepsilon}
\def\vphi{\varphi}
\def\wt{\mbox{\sl wt}\,}
\def\Z{{\bf Z}}

\title{ Energy Functions in Box Ball Systems }

\author{
Kaori Fukuda\thanks{
Department of Mathematics, Faculty of Science,
Kobe University, Rokko, Kobe 657-8501, Japan},
Masato Okado\thanks{
Department of Informatics and Mathematical Science,
Graduate School of Engineering Science,
Osaka University, Toyonaka, Osaka 560-8531, Japan}
and Yasuhiko Yamada$^*$
}

\date{}
\maketitle

\begin{abstract}
\noindent
The box ball system is studied in the crystal theory formulation.
New conserved quantities and the phase shift of the 
soliton scattering are obtained by considering 
the energy function (or $H$-function) in the combinatorial $R$-matrix.
\end{abstract}

\section{Introduction}
The box ball system (BBS) \cite{TS,T}
is one of the important classes of 
soliton cellular automata, namely
integrable discrete dynamical systems with solitons.
Recently, an interesting relation between the BBS and
the crystal theory \cite{K,KMN} has been revealed. 
In particular, it is proved in \cite{TNS}
that the scattering rule of two solitons 
in BBS is identical with the combinatorial $R$-matrix \cite{NY}.
Motivated by this and the work \cite{HIK}, 
the BBS is reformulated and generalized \cite{HKT} in terms of the 
crystal theory. In this formulation,
the time evolution is described by the row-to-row transfer 
matrix obtained as a product of the combinatorial $R$-matrices.

In the affine crystal theory, the 
combinatorial $R$-matrix consists of the `isomorphism' and
the `energy function' or $H$-function. 
In the previous works, however, the $H$-function
was not taken into account.
In this paper, we study the BBS based on the crystal 
formulation and clarify the role of the $H$-function.
The main results are as follows:
\begin{itemize}
\item{The phase shift of the soliton scattering is described by
the $H$-function. }
\item{New conserved quantities are obtained by 
including the $H$-function into the row-to-row transfer matrices.}
\end{itemize}

In the simplest BBS, that corresponds to the $n=2$ case in our paper,
it is known \cite{TTS} that a Young diagram obtained
by the Schensted bumping procedure \cite{F} is conserved
under the time evolution.
Their conserved quantities are related to the number of solitons,
and the quantities in this paper can be regarded as a generalization 
to theirs.
We note that equivalent conserved quantities are also obtained from those
of the hungry Toda molecule equation through ultra-discretization
\cite{TNS},\cite{N}.

The paper is organized as follows.
In section 2, we recapitulate necessary facts from the
crystal theory and reformulate the BBS.
In section 3, we construct conserved quantities.
The main theorem is given in section 4, where
the scattering of solitons is studied.
Section 5 is devoted for discussions.

\section{Preliminaries}
In this section we recapitulate necessary facts for the 
$U'_q(\slchap)$-crystal $B_l$ and then reformulate 
the box ball system.

\subsection{Crystal $B_l$}
Fix an integer $n\in\Z_{\ge2}$. For an arbitrary $l\in\Z_{>0}$ we set
\[
B_l=\{(\nu_1,\nu_2,\cd,\nu_l)\mid\nu_j\in\{1,2,\cd,n\},
\nu_1\le\nu_2\le\cd\le\nu_l\}.
\]
$B_l$ can be identified with the set of semi-standard tableaux of 
shape $(l)$ with 
letters in $\{1,2,\cd,n\}$. For $i=0,1,\cd,n-1$ we introduce 
operators $e_i,f_i:
B_l\longrightarrow B_l\sqcup\{0\}$ by the following rule.
For $i=1,\cd,n-1$, $e_i(\nu_1,\nu_2,\cd,\nu_l)$ is obtained by replacing
the leftmost $i+1$ in $(\nu_1,\nu_2,\cd,\nu_l)$ with $i$ and
$f_i(\nu_1,\nu_2,\cd,\nu_l)$ by replacing the rightmost $i$ with $i+1$.
If we cannot find such letter to be replaced, the result of the action
should be considered as $0$. For $i=0$ the actions are given by 
\begin{eqnarray*}
e_0(\nu_1,\nu_2,\cd,\nu_l)&=&\delta_{\nu_11}(\nu_2,\cd,\nu_l,n),\\
f_0(\nu_1,\nu_2,\cd,\nu_l)&=&\delta_{\nu_ln}(1,\nu_1,\cd,\nu_{l-1}).
\end{eqnarray*}
Notice if $f_ib=b'$ for $b,b'\in B_l$, then $b=e_ib'$.
$B_l$ is the crystal base \cite{K} of the $l$-th symmetric tensor 
representation of the quantum affine algebra $U'_q(\slchap)$.
For precise definitions, see \cite{K,KMN}.

\begin{example}
$n=4,l=5$
\begin{eqnarray*}
e_2(11334)&=&(11234)\\
f_2(11334)&=&0\\
f_0(11334)&=&(11133)
\end{eqnarray*}
\end{example}

Drawing a line from $b$ to $b'$ with color $i$ as 
$b\stackrel{i}{\longrightarrow}b'$, $B_l$ turns out a colored oriented graph
called crystal graph.

\begin{example}

\begin{itemize}
\item[(1)] n: arbitrary

$B_1:$
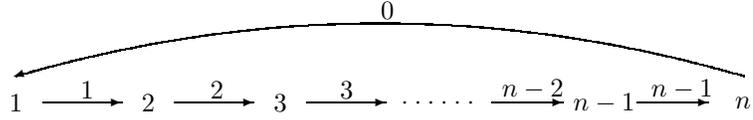
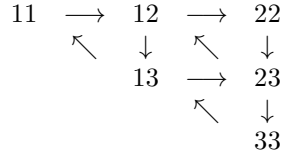
\[
\begin{picture}(300,80)(0,0)
\put(10,10){
\bezier{500}(10,20)(147.5,60)(285,20)
\put(10,20){\vector(-3,-1){1}}
\multiput(20,10)(50,0){3}{\vector(1,0){30}}
\put(5,5){\makebox(10,10){$1$}}
\put(55,5){\makebox(10,10){$2$}}
\put(105,5){\makebox(10,10){$3$}}
\put(32,10){\makebox(10,10){$1$}}
\put(81,10){\makebox(10,10){$2$}}
\put(130,10){\makebox(10,10){$3$}}
\multiput(157,10)(5,0){6}{\circle*{1}}
\multiput(190,10)(55,0){2}{\vector(1,0){27}}
\put(225,5){\makebox(15,10){$n-1$}}
\put(280,5){\makebox(10,10){$n$}}
\put(198,10){\makebox(15,10){$n-2$}}
\put(257,10){\makebox(10,10){$n-1$}}
\put(145.5,40){\makebox(10,10){$0$}}
}
\end{picture}
\]

\item[(2)] n=3

$B_2:$
\[
\begin{array}{ccccc}
11&\longrightarrow&12&\longrightarrow&22\\
  &\nwarrow&\downarrow&\nwarrow&\downarrow\\
  &&13&\longrightarrow&23\\
  &&&\nwarrow&\downarrow\\
  &&&&33
\end{array}
\]
\end{itemize}
Here the colors of the arrows $\longrightarrow,\downarrow,\nwarrow$
are $1,2,0$, respectively.
\end{example}

\subsection{Tensor product} \label{subsec:tensor}
For $b\in B_l$ we define
\[
\veps_i(b)=\max\{m\ge0\mid e_i^m b\ne0\},\quad
\vphi_i(b)=\max\{m\ge0\mid f_i^m b\ne0\}.
\]
If $b=(\nu_1,\cd,\nu_l)\in B_l$, $\veps_i(b)$ (resp. $\vphi_i(b)$)
is the number of $i+1$ (resp. $i+n\delta_{i0}$) in $(\nu_1,\cd,\nu_l)$.
For crystals $B_l$ and $B_{l'}$ the tensor product 
$B_l\ot B_{l'}=\{b\ot b'\mid b\in B_l,b'\in B_{l'}\}$ is defined.
The operators $e_i,f_i$ act on $B_l\ot B_{l'}$ as 
\begin{eqnarray}
e_i(b\ot b')&=&\left\{
\begin{array}{ll}
e_i b\ot b'&\mbox{ if }\vphi_i(b)\ge\veps_i(b')\\
b\ot e_i b'&\mbox{ if }\vphi_i(b) < \veps_i(b'),
\end{array}\right. \label{eq:tensor-e}\\
f_i(b\ot b')&=&\left\{
\begin{array}{ll}
f_i b\ot b'&\mbox{ if }\vphi_i(b) > \veps_i(b')\\
b\ot f_i b'&\mbox{ if }\vphi_i(b)\le\veps_i(b').
\end{array}\right. \label{eq:tensor-f}
\end{eqnarray}
Here $b\ot 0$ and $0\ot b'$ are understood to be $0$.
Two crystals $B_l\ot B_{l'},B_{l'}\ot B_l$ are known to be isomorphic,
{\it i.e.} the crystal graphs are the same. We denote the image of $b\ot b'$
under this isomorphism by $\btilde'\ot\btilde$. Namely, we have 
\begin{eqnarray}
B_l\ot B_{l'}&\simeq& B_{l'}\ot B_l \label{eq:iso}\\
b\ot b'&\mapsto&\btilde'\ot \btilde. \nonumber
\end{eqnarray}

We introduce so-called `signature rule' to calculate the actions of
$e_i,f_i$ on the tensor product. It is most effective on the 
multi-component tensor product. Consider an element 
$b_1\ot b_2\ot\cd\ot b_m$ of $B_{l_1}\ot B_{l_2}\ot\cd\ot B_{l_m}$.
To this element we associate an $i$-signature:
\[
\eta_i=\ubb(-,{\veps_i(b_1)})\ubb(+,{\vphi_i(b_1)})
\ubb(-,{\veps_i(b_2)})\ubb(+,{\vphi_i(b_2)})\cd
\ubb(-,{\veps_i(b_m)})\ubb(+,{\vphi_i(b_m)})
\]
Set $\eta_i^{(0)}=\eta_i$. We construct $\eta_i^{(k)}$ from $\eta_i^{(k-1)}$ 
by deleting an adjacent $(+-)$ pair in $\eta_i^{(k-1)}$ ($k=1,2,\cd$).
Then we end up with the signature of form 
$\overline{\eta}_i=\ubb(-,\alpha)\ubb(+,\beta)$ ($\alpha,\beta\ge0$).
We call such signature reduced. Of course, there are many ways of deleting
$(+-)$ pairs, but the reduced signature does not depend on them.
Now the action of $e_i$ (resp. $f_i$) on $B_{l_1}\ot B_{l_2}\ot\cd\ot
B_{l_m}$ is obtained by changing the rightmost $-$ (resp. leftmost $+$)
of $\overline{\eta}_i$ to $+$ (resp. $-$). If there is no sign to be 
changed, the action is understood as $0$. It is easy to see that in the case
of 2 components, this rule agrees with 
(\ref{eq:tensor-e}),(\ref{eq:tensor-f}).
For details, see next example.

\begin{example}
$n=4,m=3$

Consider an element $b=(1223)\ot(112)\ot(24)\in B_4\ot B_3\ot B_2$.
The $1$-signature is given as follows.
\[
\begin{array}{ccccccc}
b&=&(1223)&\ot&(112)&\ot&(24)\\
\eta_1&=&--+&&-++&&-
\end{array}
\]
The reduced signature is 
\[
\begin{array}{cccccc}
&&1&1&2\\
\overline{\eta}_1&=&-&-&+&,
\end{array}
\]
where the upper number signifies the component of the tensor product
the sign belonged to. Therefore, we have
\begin{eqnarray*}
e_1b&=e_1(1223)\ot(112)\ot(24)&=(1123)\ot(112)\ot(24),\\
f_1b&=(1223)\ot f_1(112)\ot(24)&=(1223)\ot(122)\ot(24).
\end{eqnarray*}
\end{example}

\subsection{Energy function}
Next consider a $\Z$-valued function $H$ on $B_l\ot B_{l'}$ satisfying 
the following property: For any $b\in B_l,b'\in B_{l'}$ and $i$ such that
$e_i(b\ot b')\ne0$,
\[
H(e_i(b\ot b'))=\left\{%
\begin{array}{ll}
H(b\ot b')+1&\mbox{ if }i=0,\ \vphi_0(b)\geq\veps_0(b'),\ 
\vphi_0(\btilde')\geq\veps_0(\btilde),\\
H(b\ot b')-1&\mbox{ if }i=0,\ \vphi_0(b)<\veps_0(b'),\ 
\vphi_0(\btilde')<\veps_0(\btilde),\\
H(b\ot b')&\mbox{ otherwise }.
\end{array}\right.
\]
$H$ is known to exist and unique up to additive constant.
The existences of the isomorphism (\ref{eq:iso}) and energy function $H$
are guaranteed by the existence of the $R$-matrix. See \cite{KMN}.

Explicit rules to calculate the isomorphism and energy function 
are obtained in \cite{NY}. 
Here we summarize the procedure to obtain them.
We assume $l\ge l'$.
Let $b \ot b'$ be an element in $B_l \ot B_{l'}$
such as $b=(\nu_1,\ldots,\nu_l)$ and $b'=(\nu'_1,\ldots,\nu'_{l'})$.
Let $x_i$ (resp. $y_i$) be the number of $i$ in $b$ (resp. $b'$).
We represent $b \ot b'$ by the two column diagram.
Each column has $n$ rows, enumerated as 1 to $n$ from the top to the bottom.
We put $x_i$ (resp. $y_i$) letters $i$ in the $i$-th row
of the left (resp. right) column.
The rule to obtain the energy function $H$ and the
isomorphism is as follows.
\begin{itemize}
\item[(1)]
Pick any letter, say $i$, in the right column and connect it
with a letter $j$ in the left column by a line.
The partner $j$ is chosen so that $j=\max \{\nu_k \ \vert \ \nu_k <i \}$. 
If there is no such $j$, we return to the bottom and
the partner $j$ is chosen so that $j=\max \{\nu_k \}$. 
In the latter case, we call such a pair or line ``winding".
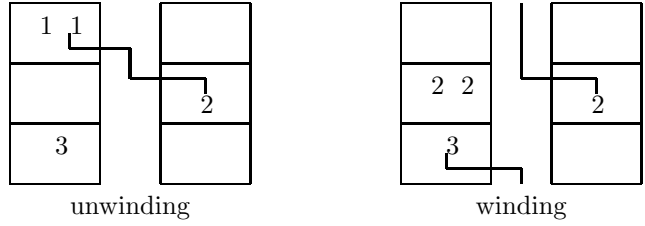
\begin{center}
\unitlength=0.8mm
\begin{picture}(105,35)(0,-5)
\multiput(0,0)(65,0){2}{
        \multiput(0,0)(25,0){2}{
                \multiput(0,0)(0,10){4}{\line(1,0){15}}
                \multiput(0,0)(15,0){2}{\line(0,1){30}}
        }
        \put(31.7,11.7){2}
}
\put( 5,25){1}
\put(10,25){1}
\put(7.5,5){3}
\put(70,15){2}
\put(75,15){2}
\put(72.5,5){3}
%
\thicklines                        
\put(20,17.5){\line(1,0){12.5}}
\put(10,22.5){\line(1,0){10}}
\put(20,22.5){\line(0,-1){5}}
\put(10,25){\line(0,-1){2.5}}
\put(32.5,15){\line(0,1){2.5}}
\put(  85,17.5){\line(1,0){12.5}}
\put(72.5, 2.5){\line(1,0){12.5}}
\put(  85, 2.5){\line(0,-1){2.5}}
\put(72.5,   5){\line(0,-1){2.5}}
\put(  85,17.5){\line(0,1){12.5}}
\put(97.5,  15){\line(0,1){2.5}}
\put(20,-2){\makebox(0,0)[t]{unwinding}}
\put(85,-2){\makebox(0,0)[t]{winding}}
\end{picture}
\end{center}
\item[(2)]
Repeat the procedure (1) for the remaining unconnected letters
$(l-1)$-times.
\item[(3)]
The isomorphism is obtained by sliding
the remaining $(k-l)$ unpaired letters in the left column to
the right.
\item[(4)]
The value of the energy function is (-1) times the number of the 
``unwinding" pairs.
\end{itemize}
Note that we normalized $H$ so that we have $H((n^l)\ot(n^{l'}))=0$.
Here and later, $(n^l)$ means $(n,n,\cd,n)\in B_l$.
When $l'=1$, the rule simplifies to the following:

If there exists $k$ such that $\nu_k<\nu'$, 
then 
\begin{eqnarray}
&&(\nu_1,\cd,\nu_l)\ot(\nu')\mapsto
(\nu_j)\ot(\nu_1,\cd,\mathop{\nu'}^j,\cd,\nu_l),
\label{eq:iso1}\\
&&H((\nu_1,\cd,\nu_l)\ot(\nu'))=-1, \nonumber
\end{eqnarray}
where $j=\max \{ k \ \vert \ \nu_k < \nu' \}$.
Otherwise,
\begin{eqnarray}
&&(\nu_1,\cd,\nu_l)\ot(\nu')\mapsto
(\nu_l)\ot(\nu',\nu_1,\cd,\nu_{l-1}),
\label{eq:iso2}\\
&&H((\nu_1,\cd,\nu_l)\ot(\nu'))=0. \nonumber
\end{eqnarray}

\begin{example}
$l=4,l'=2$
\begin{eqnarray*}
&(1123)\ot(23)\simeq(12)\ot(1233)\qquad&H=-2\\
&(1123)\ot(12)\simeq(13)\ot(1122)\qquad&H=-1\\
&(2344)\ot(12)\simeq(44)\ot(1223)\qquad&H=0
\end{eqnarray*}
\end{example}

\subsection{Yang-Baxter equation}
Let us define the affinization $\Aff(B_l)$ of the crystal $B_l$.
We introduce an indeterminate $z$ (the spectral parameter)
and set 
\[
\Aff(B_l)=\{z^d b\mid d\in\Z,b\in B_l\}.
\]
Thus $\Aff(B_l)$ is an infinite set. $z^0b\in\Aff(B_l)$ will often be 
written as $b$.

\begin{definition}
A combinatorial $R$-matrix for a crystal $B_l\ot B_{l'}$ is a map
$R:\Aff(B_l)\ot\Aff(B_{l'})\longrightarrow\Aff(B_{l'})\ot\Aff(B_l)$
given by 
\[
R(z^d b\ot z^{d'} b')=z^{d'+H(b\ot b')}\btilde'\ot z^{d-H(b\ot b')}\btilde.
\]
\end{definition}

The following result is a direct consequence of the ordinary
({\it i.e.} not combinatorial) Yang-Baxter equation.

\begin{proposition}[Yang-Baxter equation] \label{prop:YBeq}
The following equation holds on $\Aff(B_l)\ot\Aff(B_{l'})\ot\Aff(B_{l''})$.
\[
(R\ot1)(1\ot R)(R\ot1)=(1\ot R)(R\ot1)(1\ot R)
\]
\end{proposition}

\subsection{Box ball system}
Here and in what follows, we set $B=B_1$ and consider the crystal
$B^{\ot L}$ for sufficiently large $L$. The elements of $B^{\ot L}$ 
we have in mind are of the following form.
\[
\cd\ot(n)\ot(n)\ot\cd\ot\ot(n)\ot(\nu_1)\ot(\nu_2)\ot\cd\ot(\nu_k)\ot(n)
\ot(n)\ot\cd\ot(n)\ot\cd,
\]
where $\nu_1,\nu_2,\cd,\nu_k\in\{1,2,\cd,n\}$. Namely, relatively
few elements are non $(n)$, and almost all are $(n)$. In the propositions
below, we embed, if necessary, $B^{\ot L}$ into $B^{\ot L'}$ ($L<L'$) by
\begin{eqnarray*}
B^{\ot L}&\hookrightarrow& B^{\ot L'}\\
b_1\ot\cd\ot b_L&\mapsto& b_1\ot\cd\ot b_L\ot
\underbrace{(n)\ot\cd\ot(n)}_{L'-L}.
\end{eqnarray*}

\begin{lemma} \label{lem:1}
By iterating the isomorphism $B_l\ot B\simeq B\ot B_l$, we consider a map
\begin{eqnarray*}
B_l\ot B\ot\cd\ot B&\mathop{\longrightarrow}^\sim&
B\ot\cd\ot B\ot B_l\\
(n^l)\ot b_1\ot\cd\ot b_L&\mapsto&\btilde_1\ot\cd\ot\btilde_L\ot\btilde.
\end{eqnarray*}
Then there exists an integer $L_0$ such that $\btilde=(n^l)$ for $L\ge L_0$.
\end{lemma}
Taking sufficiently large $L$ such that the above lemma holds, we define
a map $T_l:B^{\ot L}\longrightarrow B^{\ot L}$ by 
$b_1\ot\cd\ot b_L\mapsto\btilde_1\ot\cd\ot\btilde_L$.

\begin{lemma} \label{lem:2}
For a fixed element of $B^{\ot L}$, there exists an integer $l_0$ such that
$T_l=T_{l_0}$ for any $l\ge l_0$.
\end{lemma}
Both lemmas are obvious from (\ref{eq:iso1}), (\ref{eq:iso2}).

\vskip5mm
Following \cite{HKT}, we reformulate the box ball system. We consider
such an element of $B^{\ot L}$, often called a `state', as described in 
the beginning of this subsection. Lemma \ref{lem:1} and \ref{lem:2} enables
us to define an operator $T=\lim_{l\rightarrow\infty}T_l$ on the space of 
states.
Application of $T$ gives rise to a transition of state. Thus it can be 
regarded as describing a certain dynamical system, in which $T$ plays
the role of `time evolution'. 
{}From the same reason, $T_l$ may also be viewed as a time evolution.
(In this paper, time evolution means the one by $T$ unless otherwise stated.)

Under the condition assumed in Lemma \ref{lem:1},
the sequence  $b_1\ot\cd\ot b_L$ is determined uniquely from
$\btilde_1\ot\cd\ot\btilde_L$ and $\btilde=(n^l)$, since the
isomorphisms $B_l \ot B \simeq B \ot B_l$ are bijective.
Hence, the time evolutions $T_l$ ($1 \geq l$) are invertible.
The inverse map of the isomorphism (\ref{eq:iso1}),(\ref{eq:iso2}) is
given explicitly by the following rule.

If there exists $k$ such that $\nu_k>\nu'$, 
then 
\begin{equation}
(\nu')\ot(\nu_1,\cd,\nu_l)\mapsto
(\nu_1,\cd,\mathop{\nu'}^j,\cd,\nu_l)\ot(\nu_j),
\end{equation}
where $j=\min \{ k \ \vert \ \nu_k > \nu' \}$.
Otherwise,
\begin{equation}
(\nu')\ot(\nu_1,\cd,\nu_l)\mapsto
(\nu_2,\cd,\nu_{l},\nu')\ot(\nu_1).
\end{equation}

As we see in the next section, these time evolutions
$T_l$ ($l \geq 1$) form a commuting family of operators
and play an important role in our paper.

\section{Conservation laws}

In this section we construct conserved quantities under the time evolutions.

\subsection{Definition of $E_l$}
Fix sufficiently large $L$ and consider a composition of the 
combinatorial $R$-matrices
\[
\R_l=R_{L\,L+1}\cd R_{23}R_{12}:
\Aff(B_l)\ot\Aff(B)^{\ot L}\longrightarrow\Aff(B)^{\ot L}\ot\Aff(B_l).
\]
Here $R_{i\,i+1}$ signifies that the $R$-matrix acts on the $i$-th and
$(i+1)$-th components of the tensor product.
Applying $\R_l$ to an element $(n^l)\ot p$ ($p=b_1\ot\cd\ot b_L$), we
have 
\begin{eqnarray*}
\R_l((n^l)\ot p)
&=&
z^{H_1}\btilde_1\ot 
z^{H_2}\btilde_2\ot\cd\ot
z^{H_L}\btilde_L\ot z^{E_l(p)}(n^l),\\
E_l(p)&=&-\sum_{j=1}^L H_j, \quad
H_j=H(b^{(j-1)}\ot b_j),
\end{eqnarray*}
where $b^{(j)}$ ($0\le j<L$) is defined by
\begin{eqnarray*}
B_l\ot\underbrace{B\ot\cd\ot B}_j&\simeq&\underbrace{B\ot\cd\ot B}_j\ot B_l\\
(n^l)\ot b_1\ot\cd\ot b_j&\mapsto&
\btilde_1\ot\cd\ot\btilde_j\ot b^{(j)}.
\end{eqnarray*}

\subsection{Conservation of $E_l$}

\begin{proposition}
As maps from $\Aff(B_l)\ot\Aff(B_{l'})\ot\Aff(B)^{\ot L}$
to $\Aff(B)^{\ot L}\ot\Aff(B_{l'})\ot\Aff(B_l)$, the following are
identical. 
\begin{equation}
\R_{l'}\R_lR_{12}=R_{L+1\,L+2}\R_l\R_{l'} \label{eq:commuting}
\end{equation}
Here $\R_l$ (resp. $\R_{l'}$) acts on $\Aff(B_l)\ot\Aff(B)^{\ot L}$
(resp. $\Aff(B_{l'})\ot\Aff(B)^{\ot L}$), and acts on the other component
as identity.
\end{proposition}

\Proof
Use the Yang-Baxter equation repeatedly.
\qed

\begin{theorem}
For an element $p\in B^{\ot L}$, we have
\begin{itemize}
\item[(1)] $T_lT_{l'}(p)=T_{l'}T_l(p)$.
\item[(2)] $E_l(T_{l'}(p))=E_l(p)$. In particular, $E_l(T(p))=E_l(p)$.
\end{itemize}
\end{theorem}

\Proof
Consider an element $(n^l)\ot(n^{l'})\ot p$ ($p=b_1\ot\cd\ot b_L$) of
$\Aff(B_l)\ot\Aff(B_{l'})\ot\Aff(B)^{\ot L}$.
Apply both sides of (\ref{eq:commuting}) to it. Then one gets
\begin{eqnarray*}
LHS&=&z^{H_1+\tilde{H}_1}\tilde{\btilde}_1\ot\cd\ot 
z^{H_L+\tilde{H}_L}\tilde{\btilde}_L
\ot z^{E_{l'}(\tilde{p})}(n^{l'})\ot z^{E_l(p)}(n^l),\\
RHS&=&z^{H'_1+\hat{H}_1}\hat{\hat{b}}_1\ot\cd\ot 
z^{H'_L+\hat{H}_L}\hat{\hat{b}}_L
\ot z^{E_{l'}(p)}(n^{l'})\ot z^{E_l(\hat{p})}(n^l).
\end{eqnarray*}
Here $\tilde{\btilde}_1\ot\cd\ot\tilde{\btilde}_L=T_{l'}T_l(p),
\hat{\hat{b}}_1\ot\cd\ot\hat{\hat{b}}_L=T_lT_{l'}(p)$ and 
$\tilde{p}=T_l(p),\hat{p}=T_{l'}(p)$. 
$H_j,\tilde{H}_j,H'_j,\hat{H}_j$ ($1\le j \le L$) are values of the 
energy functions, which are not necessary here.
(\ref{eq:commuting}) obviously ensures (1) and (2).
\qed

\section{Solitons}

\subsection{Soliton}
A state of the following form is called an $m$-soliton state of
length $l_1,l_2,\ldots,l_m$,
\begin{equation}\label{eq:soliton-state}
...[l_1]........[l_2]..... \cdots .....[l_m]....... .
\end{equation}
Here $..[l]..$ denotes a local configuration such as
\[
\cd\ot n \ot n \ot
\underbrace{\nu_1\ot\nu_2\ot\cd\ot\nu_l}_{l}\ot n \ot n \ot\cd,
\quad
(n>\nu_1 \geq \nu_2 \geq \cdots \geq \nu_l\geq 1),
\]
sandwiched by sufficiently many $n$'s.
(In this section, the parentheses in the state are omitted.)

\begin{lemma} \label{lem:single-soliton}
Let $p$ be a 1-soliton state of length $l$, then
\begin{itemize}
\item[(1)]
The $k$-th conserved quantity of $p$ is given by $E_k(p)=\min(k,l)$.
\item[(2)]
The state $T_k(p)$ is obtained by the rightward shift by
$E_k(p)$ lattice steps.
\end{itemize}
\end{lemma}
\Proof
(1) Recall that the conserved quantity $E_k$ is a sum of local 
$H$-functions $-H_j=-H(b^{(j-1)} \ot b_j) \in \{0,1\}$.
By the rule (2.2) and (2.3), one has $-H_j=1$ if and only if
$j=J+1,J+2,\ldots,J+\min(k,l)$, where $J$ is the position of 
the last component $\nu_l$ in the tensor product $B^{L}$.
Hence $E_k=\min(k,l)$.
Similarly, the statement (2) follows from the rule (2.2) and (2.3).
\qed

\noindent
\begin{example}
The time evolution\footnote{
The convention of the box ball game here is different from
\cite{TS,T} in that, bigger numbers are transferred before
smaller ones.}
of 1 soliton of length $3$.
Here and in what follows $.$ means $n$ ($n>3$).
{\footnotesize{
\[
\matrix{
.&.&.&.&3&3&2&.&.&.&.&.&.&.&.&.&.&.&.&.\cr
.&.&.&.&.&.&.&3&3&2&.&.&.&.&.&.&.&.&.&.\cr
.&.&.&.&.&.&.&.&.&.&3&3&2&.&.&.&.&.&.&.\cr
.&.&.&.&.&.&.&.&.&.&.&.&.&3&3&2&.&.&.&.
}
\]
}}
\end{example}

\begin{definition}
For any state $p$, the numbers $N_l=N_l(p)$ ($l=1,2,\ldots$)
are defined by 
\begin{eqnarray*}
&&E_l=\sum_{k \geq 1} \min(k,l) N_k, \quad E_0=0, \\
&&N_l=-E_{l-1}+2 E_{l}-E_{l+1}.
\end{eqnarray*}
\end{definition}
By Lemma \ref{lem:single-soliton}, we have
\begin{proposition}
For $m$-soliton state (\ref{eq:soliton-state}),
$N_l$ is the number of solitons of length $l$,
\[
N_l=\sharp \{ j \ \vert \  l_j=l \}.
\]
\end{proposition}
This proposition implies the stability of solitons,
since the numbers $E_l(p)$, and hence $N_l(p)$, are conserved.

\subsection{Scattering of solitons}

The following is an example of the time evolution 
(for $t=0,1,\ldots,6$) of a state
which shows the scattering of three solitons of length
3,2 and 1.
\begin{example}\label{ex:3soliton-scattering}
{\footnotesize{
\[
\matrix{
3&3&2&.&.&.&1&1&.&.&.&2&.&.&.&.&.&.&.&.&.&.&.&.&.&.\cr
.&.&.&3&3&2&.&.&1&1&.&.&2&.&.&.&.&.&.&.&.&.&.&.&.&.\cr
.&.&.&.&.&.&3&3&.&.&2&1&1&2&.&.&.&.&.&.&.&.&.&.&.&.\cr
.&.&.&.&.&.&.&.&3&3&.&.&.&1&2&2&1&.&.&.&.&.&.&.&.&.\cr
.&.&.&.&.&.&.&.&.&.&3&3&.&.&1&.&.&2&2&1&.&.&.&.&.&.\cr
.&.&.&.&.&.&.&.&.&.&.&.&3&3&.&1&.&.&.&.&2&2&1&.&.&.\cr
.&.&.&.&.&.&.&.&.&.&.&.&.&.&3&.&3&1&.&.&.&.&.&2&2&1
} 
\]
}}
\end{example}
We introduce a labeling of solitons of length $l$ using
$\Aff(B_l)$ for the lower rank algebra $U'_q(\widehat{\goth{sl}}_{n-1})$.
(See also \cite{HKT} for such an identification.) 
Suppose there is a soliton of length $l$
$\cd \ot \nu_1 \ot \nu_2 \ot \cdots \ot \nu_l \ot \cd$
($n>\nu_1 \geq \nu_2 \geq \cdots \geq \nu_l\geq 1$) at time $t$. 
Say it is at position $\gamma(t)$, if 
$(\nu_1)$ is in the $\gamma(t)$-th tensor component of $B^{\ot L}$. 
{}From Lemma \ref{lem:single-soliton} (2), the position $\gamma(t)$ 
under the time evolution $T_k$ is
given by $\gamma(t)=\min(k,l)t+\gamma$ ($\min(k,l)$ is the velocity and 
$\gamma$ is the phase)
unless it interacts with other solitons. 
To such a soliton we associate an
element $z^{-\gamma}(\nu_l,\nu_{l-1},\cd,\nu_1)\in\Aff(B_l)$ for 
$U'_q(\widehat{\goth{sl}}_{n-1})$.

Now consider a state of $m$ solitons illustrated as below.
\[
...[l_1].....[l_2]...............[l_m]................
\]
We assume solitons are separated enough from each other and
$l_1>l_2>\cd>l_m$. Since longer solitons move faster,
we can expect that the state turns out to be 
\[
.......[l_m]...........[l_2].....[l_1]................
\]
after sufficiently many time evolutions. 
(The proof of this fact is found in that of the main theorem.)
We represent such a scattering process as
\begin{eqnarray*}
&&z^{c_1}b_1\ot z^{c_2}b_2\ot\cd\ot z^{c_m}b_m \\
&\mapsto&z^{c'_m}b'_m\ot\cd\ot z^{c'_2}b'_2\ot z^{c'_1}b'_1
\end{eqnarray*}
Here $z^{c_j}b_j,z^{c'_j}b'_j$ are elements of $\Aff(B_{l_j})$ 
under the identification in the previous paragraph.
With these notations, the scattering process 
in Example \ref{ex:3soliton-scattering} is described as
\begin{equation} \label{eq:scattering-example}
z^{0}(2,3,3)\ot z^{-6}(1,1)\ot z^{-11}(2)
\mapsto
z^{-8}(3)\ot z^{-4}(1,3)\ot z^{-5}(1,2,2).
\end{equation}

Let us recall some useful fact derived from representation theory.
Note that $U'_q(\slchap)$ contains $U_q(\goth{sl}_{n-1})$
as subalgebra. This fact can be translated into the language of crystals
and guarantees that 
\[
T_k\mbox{ commutes with }e_i,f_i\;(i=1,2,\cd,n-2)\mbox{ on }
\Aff(B_{l_1})\ot\cd\ot\Aff(B_{l_m}).
\]
Here the actions of $e_i,f_i$ ($i=1,2,\cd,n-2$) on the multi-component
tensor product can be calculated using signature rule explained 
in section \ref{subsec:tensor}. By the actions, the power of $z$ in an element
of $\Aff(B_{l_j})$ is unaffected. We call this property 
$U_q(\goth{sl}_{n-1})$-invariance. 
This property is also used to prove our theorem.
For instance, if we admit 
(\ref{eq:scattering-example}), we have
\begin{eqnarray*}
&&e_2(z^{0}(2,3,3)\ot z^{-6}(1,1)\ot z^{-11}(2))=
z^{0}(2,2,3)\ot z^{-6}(1,1)\ot z^{-11}(2)\\
&\mapsto&
e_2(z^{-8}(3)\ot z^{-4}(1,3)\ot z^{-5}(1,2,2))=
z^{-8}(3)\ot z^{-4}(1,2)\ot z^{-5}(1,2,2).
\end{eqnarray*}

Let us explain the $U_q(\goth{sl}_{n-1})$-invariance in terms of 
tableau combinatorics. Consider a state, read letters from right to left
and remove all $n$'s. Suppose the resulting sequence of letters is 
$w=i_1\cd i_s$ ($1\le i_k\le n-1$). From such a sequence (or word)
one can define a semi-standard tableau
$(\cd((i_1\leftarrow i_2)\leftarrow i_3)\cd\leftarrow i_s)$ by the
`row bumping' algorithm \cite{F}. The $U_q(\goth{sl}_{n-1})$-invariance 
assures that the resulting tableau is invariant under the time evolution.
In the example (\ref{eq:scattering-example}), the words before and after
the scattering are given by $w=211233$ and $w'=122133$, and the corresponding
tableaux are the same:
\[
\begin{array}{ccccc}
1&1&2&3&3\\
2\\
\end{array}
\]
Of course, this invariance is valid in the intermediate stage of scattering.

The main results of this paper are the following.

\begin{theorem}\label{th:main}
\begin{itemize}
\item[(1)] The two body scattering of solitons of length $l_1$ and
length $l_2$ ($l_1 > l_2$) under the time evolution $T_r$ ($r>l_2$)
is described by the combinatorial $R$-matrix:
\begin{eqnarray*}
\Aff(B_{l_1})\ot \Aff(B_{l_2})&\simeq&\Aff(B_{l_2})\ot \Aff(B_{l_1})\\
z^{c_1}b_1 \ot z^{c_2}b_2 &\mapsto& 
z^{c_2+\delta}\btilde_2 \ot z^{c_1-\delta}\btilde_1,
\end{eqnarray*}
with $\delta=2 l_2+H(b_1 \ot b_2)$.

\item[(2)] The scattering of solitons is factorized into
two body scatterings.
\end{itemize}
\end{theorem}

Some remarks may be in order. 

\noindent
1) The combinatorial $R$-matrix in Theorem
\ref{th:main} has an extra term $2l_2$ in the power of $z$.
However, the Yang-Baxter equation (Proposition \ref{prop:YBeq})
holds as it is.

\noindent
2) Although we do not consider the cases where there exist
solitons with same length, some part of the results can be
generalized to such situations. 
Let us consider a state 
$$
..\underbrace{[k]....[k]}_{N_k}..... \cdots
..\underbrace{[2]....[2]}_{N_2}.....
..\underbrace{[1]....[1]}_{N_1}..... ,
$$
consisting of $N_j$ solitons of length $j$ ($1 \leq j \leq k$).
By the same argument given in the following proof, we can show that
the scattering of these solitons is factorized into 
scattering of two bunches ($1 \leq j < i \leq k$)
$$
..\underbrace{[i]....[i]}_{N_{i}}..... \cdots
..\underbrace{[j]....[j]}_{N_{j}}..... .
$$
We conjecture this scattering is described by the product of 
$N_i N_j$ combinatorial $R$ matrices
$\Aff(B_{i})\ot \Aff(B_{j})\simeq\Aff(B_{j})\ot \Aff(B_{i})$.
This conjecture is trivially true if each solitons are separated
enough.

\begin{example}
Using this theorem, the scattering process 
(\ref{eq:scattering-example}) is calculated as
\begin{eqnarray*}
&&z^{0}(2,3,3) \ot z^{-6}(1,1) \ot z^{-11}(2) 
=\left\{z^{0}(2,3,3) \ot z^{-6}(1,1) \right\} \ot z^{-11}(2) \\
&\mapsto& \left\{z^{-6+4}(3,3) \ot z^{0-4}(1,1,2) \right\} \ot z^{-11}(2)
= z^{-2}(3,3) \ot \left\{ z^{-4}(1,1,2) \ot z^{-11}(2) \right\} \\
&\mapsto& z^{-2}(3,3) \ot \left\{ z^{-11+1}(1) \ot z^{-4-1}(1,2,2) \right\}
= \left\{ z^{-2}(3,3) \ot z^{-10}(1) \right\} \ot z^{-5}(1,2,2) \\
&\mapsto& \left\{ z^{-10+2}(3) \ot z^{-2-2}(1,3) \right\} \ot z^{-5}(1,2,2) 
=z^{-8}(3) \ot z^{-4}(1,3) \ot z^{-5}(1,2,2).
\end{eqnarray*}
The result is independent of the order of the factorization
due to the Yang-Baxter equation.
\end{example}

\noindent
{\it Proof of Theorem \ref{th:main}.}
(1) Two soliton scattering rule: 
Due to the $U_q(\sln-)$-invariance, it is enough to 
check the rule for the highest weight elements
$(1^i) \ot (1^l,2^k) \in B_i \ot B_{k+l}$ ($i > l+k$), {\it i.e.} 
elements such that $e_j(b_1\ot b_2)=0$ for any $j=1,\cd,n-2$.
The corresponding state is given by
\[
\cd\ot\ub(1,i)\ot\ub(n,j)\ot\ub(2,k)\ot\ub(1,l)\ot\cd.
\]
Now consider the time evolution by $T_{k+l+1}$.
If $j>l$, $(1^i)$ moves with velocity $k+l+1$ and $(1^l,2^k)$ with
$k+l$. At some time, we arrive at the state
\[
\cd\ot\ub(1,i)\ot\ub(n,l)\ot\ub(2,k)\ot\ub(1,l)\ot\cd.
\]
By the rule (\ref{eq:iso1})(\ref{eq:iso2}), it is easy to see that
after $t$ time units $(t\le i-k-l)$ from this moment, the state becomes
\[
\cd\ot\ub(n,{(k+l+1)t})\ot\ub(1,i-t)\ot\ub(n,l)\ot\ub(2,k)\ot\ub(1,l+t)\ot\cd.
\]
After that $(t>i-k-l)$, the solitons never interact again.

We have investigated the time evolution by $T_{k+l+1}$.
However, we have the equality 
\[
T_r^b T_{k+l+1}^a = T_{k+l+1}^a T_r^b\qquad(r>k+l).
\]
If $a$ and $b$ are sufficiently large, we can reduce the observation
of the scattering by $T_r$ in the right hand side to that by $T_{k+l+1}$
in the left hand side, which we have just finished.
It is immediate to check that 
$(1^i)\ot(1^l,2^k)\mapsto(1^{l+k})\ot(1^{i-k},2^k)$ under the isomorphism
$B_i\ot B_{l+k}\simeq B_{l+k}\ot B_i$. Let $\delta_1$ (resp. $\delta_2$) be 
the phase shift of the soliton of length $i$ 
(resp. soliton of length $(l+k)$). They are given by 
\begin{eqnarray*}
\delta_1&=&(k+l+1)(i-k-l)+(k+l)+l-(k+l+1)(i-k-l)=2l+k,\\
\delta_2&=&(k+l+1)(i-k-l)-i-l-(k+l)(i-k-l)=-(2l+k).
\end{eqnarray*}
Since $H((1^i)\ot(1^l,2^k))=-k$, we have the desired result.


(2) Factorization property:
We illustrate the proof by a scattering of three solitons
of length 3,2 and 1. Other cases are similar.
Now the initial state is
\[
p=...[3]..[2]........[1]................
\]
We use a similar technique to (1).
By applying the operator $T_2^a$ ($a>>0$) 
we have scatterings of $[2]\ot[1]$ and $[3]\ot[1]$.
\[
T_2^a(p)=.........[1].........[3]..[2]...........
\]
Then by $T_3^b$ ($b>>0$) we have the 
third scattering of $[3]\ot[2]$.
\[
T_3^b T_2^a(p)=.........[1].........[2]......[3]...........
\]
After this the time evolution $T^c$ acts trivially.
On the other hand, these evolution can be rewritten as
\[
T^c T_3^b T_2^a = T_3^b T_2^a T^c.
\]
In the right hand side, we first observe three soliton scattering
caused by $T^c$ ($c>>0$) and then $T_3^b T_2^a$ act trivially.
Hence, the three soliton scattering in the right hand side
is factorized into 3 $\times$ (two soliton scattering) 
in the left hand side.
\qed

\section{Discussion}
In this paper we have studied the box ball system in the
crystal base formulation and clarified the role of $H$-functions.
The $H$-function plays two important roles in the system:
(1) to construct the conserved quantities
and (2) to describe the phase shift of the scattering.
Similar results are expected for the generalized models
defined in \cite{HKT}.

As far as we know, there is no direct proof of the equivalence between 
the original box ball systems and their crystal theoretic reformulation.
This problem is under investigation.

Another important feature of the box ball systems
is their relation to nonlinear integrable equations through the
`ultra discretization' \cite{MSTTT,TTMS,TNS}.
It is an interesting problem to understand the direct relation
between such nonlinear integrable equations and crystal theory.

\vskip0.5cm\noindent
{\bf Acknowledgements}.
The authors would like to thank G.~Hatayama, A.~Kuniba, M.~Noumi, T.~Takagi, 
T.~Takebe and T.~Tokihiro for valuable discussions.
They also thank G.~Hatayama, A.~Kuniba and T.~Takagi for sending \cite{HKT}
prior to publication.

\end{document}